\documentclass[a4paper,12pt]{article}
\usepackage{amssymb,amsmath,graphics,pstricks}

\newenvironment{dwd}{\par\noindent{\bf Proof.}}{\par\rightline{$\blacksquare$}}
\newenvironment{dww}{\par\noindent{\bf Proof of Theorem 1.}}{\par\rightline{$\blacksquare$}}
\newenvironment{dw1}{\par\noindent{\bf Proof of the left-hand side of (\ref{kwiat1}).}}{\par\rightline{$\blacksquare$}}
\newenvironment{dw2}{\par\noindent{\bf Proof of the right-hand side of (\ref{kwiat1}).}}{\par\rightline{$\blacksquare$}}
\newenvironment{dw3}{\par\noindent{\bf Proof of Theorem 2 in the case of $\eta'(0)<\infty$.}}{\par\rightline{$\blacksquare$}}
\newtheorem{theo}{Theorem}
\newtheorem{lema}{Lemma}
\newtheorem{prop}{Proposition}
\newtheorem{defi}{Definition}

\newtheorem{rema}{Remark}
\newtheorem{coro}{Corollary}

\def\be#1\ee{\begin{equation}#1\end{equation}}
\newcommand{\ba}{\begin{eqnarray} }
\newcommand{\ea}{\end{eqnarray} }
\def\bt#1\et{\begin{theo}#1\end{theo}}
\def\bl#1\el{\begin{lema}#1\end{lema}}
\def\bp#1\ep{\begin{prop}#1\end{prop}}
\def\bd#1\ed{\begin{defi}#1\end{defi}}

\def\ccB{{\cal B}}

\def\ccF{{\cal F}}

\def\ccM{{\cal M}}

\def\ccP{{\cal P}}

\def\ccS{{\cal S}}

\def\va{\varepsilon}
\def\ra{\rightarrow}

\def\pint{-\hspace{-11pt}\int}

\def\E{\mathbf{E}}
\def\P{\mathbf{P}}

\def\R{{\mathbb R}}

\def\ls{\leqslant}
\def\gs{\geqslant}
\def\for{\mbox{for}}
\begin{document}
\title {\bf On a type Sobolev inequality and its applications.}
\author{Witold Bednorz\\ \emph{Warsaw University}}
\maketitle

\begin{abstract}
Assume $\|\cdot\|$ is a norm on $\R^n$ and $\|\cdot\|_{\ast}$ its dual.
In this paper we consider the closed ball $T:=B_{\|\cdot\|}(0,r)$, $r>0$.
Suppose $\varphi$ is an Orlicz function and $\psi$ its conjugate,
we prove that for arbitrary $A,B>0$ and for each Lipschitz function $f$
on $T$ the following inequality holds
\ba
&&\sup_{s,t\in T}|f(s)-f(t)|\ls
6AB(\int^r_{0}\psi(\frac{1}{A\va^{n-1}})\va^{n-1}d\va+\nonumber\\
&&+\frac{1}{n|B_{\|\cdot\|}(0,1)|}\int_{T}\varphi(\frac{1}{B}\|\nabla f(u)\|_{\ast})du)\nonumber,
\ea
where $|\cdot|$ is the standard Lebesgue measure on $\R^n$.
This is a strengthening of the Sobolev inequality obtained in the proof of Theorem 5.1 by
M. Talagrand \cite{Tal0}. We use the inequality to state for a given concave, strictly increasing function
$\eta:\R_{+} \ra \R$, with $\eta(0) = 0$, the necessary and sufficient
condition on $\varphi$ so that each separable process $X(t)$, $t\in T$ which satisfies
$$
\|X(s)-X(t)\|_{\varphi}\ls \eta(\|s-t\|),\;\;\for\; s,t\in T
$$
is a.s. sample bounded.
\end{abstract}
\renewcommand{\thefootnote}{\noindent}
\footnote{\noindent\hspace{-.6cm}{\bf Subject classification:}  60G17, 28A99.
Partially supported by the Funds of Grant MENiN 1 P03A 01229}

\section{Introduction}

Let $\|\cdot\|$ be a norm on $\R^n$. We denote by $B_{\|\cdot\|}(x,r)$ the
closed ball with the center at $x$ and the radius $r$ with respect to 
the metric given by $\|\cdot\|$, i.e.  
$$
B_{\|\cdot\|}(x,r):=\{y\in \R^n:\;\|x-y\|\ls r\}.
$$
Let $\langle \cdot,\cdot \rangle$ be the canonical scalar product
(that is $\langle u,v \rangle:=\sum^n_{i=1}u_iv_i$, for $u,v\in \R^n$)
and $\|\cdot\|_{\ast}$ the dual norm, i.e.
$$
\|v\|_{\ast}:=\sup_{u \in B_{\|\cdot\|}(0,1)}|\langle u,v\rangle|,\;\;\for\;v\in \R^n.
$$
In this paper we consider the closed ball $T:=B_{\|\cdot\|}(0,r)$, $r>0$.
\smallskip

\noindent
We say that $\varphi:\R_{+}\ra \R$ is an Orlicz function if it is convex, strictly increasing,
$\varphi(0)=0$ and also $\lim_{x\ra 0}\varphi(x)/x=0$, $\lim_{x\ra \infty}\varphi(x)/x=\infty$.
For each Orlicz function $\varphi$ we define its conjugate
$$
\psi(x):=\sup_{y\gs 0}(xy-\varphi(y)),\;\;\for\;x\gs 0.
$$
This $\psi$ is also an Orlicz function.
Moreover, it is well known that $\varphi$ is the conjugate function for $\psi$, namely
$\varphi(x)=\sup_{y\gs 0}(xy-\psi(y))$.
The definition implies the Young inequality
\be\label{nie3}
xy\ls \varphi(x)+\psi(y),\;\;\for\; x,y\gs 0.
\ee
From now on we assume that $\varphi,\psi$ are conjugate Orlicz functions.
\smallskip

\noindent
In the paper we prove the following
Sobolev type  inequality and give its applications to the theory of
stochastic processes.
\bt\label{natul1}
For each $A,B>0$ and for each Lipschitz function $f$ on $T$
the following inequality holds
\ba
&&\sup_{s,t\in T}|f(s)-f(t)| \ls
6AB(\int^r_{0}\psi(\frac{1}{A\va^{n-1}})\va^{n-1}d\va+\nonumber\\
&&+\frac{1}{n|B_{\|\cdot\|}(0,1)|}\int_{T}\varphi(\frac{1}{B}\|\nabla f(u)\|_{\ast})du)\nonumber,
\ea
where $|\cdot|$ is the standard Lebesgue measure on $\R^n$.
\et
The above inequality is a generalization of Talagrand's result,
who obtained such inequality in the proof of Theorem
5.1, \cite{Tal0} when $\|s-t\|=\sup^n_{i=1}|s_i-t_i|$.
\smallskip

\noindent
Since for each $s,t\in\R^n$ we have $s\in B_{\|\cdot\|}(t,\|s-t\|)$, the above theorem
implies some regularity on $f$. Namely, for $\varphi(x)\equiv x^p/p$, $p>n$ we obtain the
following classical result (which can be deduced from Lemma 7.16 in \cite{Gal-Tru} 
by using Holder inequality).
\begin{coro}
Suppose $p>n$, then for each Lipschitz function $f$ on $\R^n$ the following inequality holds
$$
\sup_{s,t\in \R^n}\frac{|f(s)-f(t)|}{\|s-t\|^{1-\frac{n}{p}}}\ls
\frac{6(\frac{p-1}{p-n})^{1-\frac{1}{p}}}{(n|B_{\|\cdot\|}(0,1)|)^{1/p}} (\int_{\R^n}
\|\nabla f(u)\|_{\ast}^{p} du)^{1/p}.
$$
\end{coro}
\begin{dwd}
The conjugate function for $\varphi(x)\equiv x^p/p$ is $\psi(x)\equiv  x^{q}/q$, where $1/p+1/q=1$. 
Due to Theorem \ref{natul1} for each $A,B>0$ we have
\ba
&&|f(s)-f(t)|\ls
6AB(\int^{\|s-t\|}_{0} q^{-1}(\frac{1}{A\va^{n-1}})^q \va^{n-1}d\va+\nonumber\\
&&+\frac{1}{n|B_{\|\cdot\|}(0,1)|}\int_{B_{\|\cdot\|}(t,\|s-t\|)} p^{-1}(\frac{1}{B}\|\nabla f(u)\|_{\ast})^p du)\nonumber,
\ea
We take 
$$
A=(\int^{\|s-t\|}_{0} \va^{(1-q)(n-1)}d\va)^{1/q},\;\;
B=(\frac{1}{n|B_{\|\cdot\|}(0,1)|}\int_{\R^n} (\|\nabla f(u)\|_{\ast})^p du)^{1/p}.
$$
Consequently ($T=B_{\|\cdot\|}(0,r)$)
$$
|f(s)-f(t)|\ls 
6AB=\frac{6(\frac{p-1}{p-n})^{1/q} \|s-t\|^{1-\frac{n}{p}}}{(n|B_{\|\cdot\|}(0,1)|)^{1/p}}
(\int_{T}\|\nabla f(u) \|^p_{\ast} du)^{1/p}.
$$
It completes the proof.
\end{dwd}
Another immediate consequence of Theorem \ref{natul1} is
the sufficient condition for embedding of the Sobolev space $W^{1,\varphi}_0(T)$ into $L_{\infty}(T)$.
\begin{coro}
If for some $A>0$ we have
$\int^r_{0}\psi(\frac{1}{A\va^{n-1}})\va^{n-1}d\va<\infty$
then the space $W^{1,\varphi}_0(T)$ embeds into $L_{\infty}(T)$.
\end{coro}
The result can be deduced
from the part II of Theorem 1.1. in the paper by A. Cianchi \cite{Cia}.
\smallskip

\noindent
To explain applications to
stochastic processes we need some definitions.
Let $(K,d)$ be a compact metric space.
Denote by $\mathfrak{B}(K)$ the space of Borel bounded functions on $K$,
by $C(K)$ the Banach space of continuous functions on $K$, with sup-norm and
by $\mathrm{Lip}(K)$ the space of Lipschitz functions on $T$
with the norm
$$
\|f\|_{\mathrm{Lip}}:=\sup_{s\neq t}\frac{|f(s)-f(t)|}{d(s,t)}+D(K)^{-1}\|f\|_{\infty},
$$
where $D(K):=\sup\{d(s,t):\;s,t\in K\}$ (the diameter of $K$).
Let $\ccP(K)$ be the set of all Borel probability measures on $K$.
For each $\nu\in \ccP(K)$, $f\in \mathfrak{B}(K)$ and  $A\in\ccB(K)$ (with $\nu(A)>0$) we denote
$$
\pint_{A}f(u)\nu(du):=\frac{1}{\nu(A)}\int_A f(u)\nu(du).
$$
Suppose $X$ is a random variable, we define the Luxemburg norm $\|X\|_{\varphi}:=
\inf\{c>0:\;\E\varphi(\frac{|X|}{c})\ls 1\}$. For a fixed probability space the Banach
space $L_{\varphi}$ consists of all random variables for which $\|X\|_{\varphi}<\infty$.
\smallskip

\noindent
In this paper we consider only separable processes (for the definition
see Introduction in the book by Ledoux-Talagrand \cite{Le-Ta}).
For each separable $X(t)$, $t\in {K}$ we have the following equality
\be\label{hlyt}
\E\sup_{s,t\in K}|X(s)-X(t)|=\sup_{F\subset {K}}\E\sup_{s,t\in F}|X(s)-X(t)|,
\ee
where the supremum is taken over all finite subsets of $K$.
Let us impose the Lipschitz condition on
increments of our processes, that is for each $X(t)$, $t\in K$
we assume that
\be\label{war2}
\sup_{s,t\in K}\E\varphi(\frac{|X(s)-X(t)|}{d(s,t)})\ls 1.
\ee
This condition can be rewritten in terms of Luxemburg norms
$$
\|X(s)-X(t)\|_{\varphi}\ls d(s,t),\;\;\for\;s,t\in K.
$$
In the theory of stochastic processes a lot of effort has been put in
finding criteria for boundedness or continuity of stochastic processes. In
most of the cases they are of the following Kolmogorov type: some
assumptions on the Orlicz function $\varphi$ and the metric space $(K,d)$ are
given so that for each separable process $X$ on $K$ the condition
(\ref{war2}) implies that $X$ is bounded a.s.
\smallskip

\noindent
It is not difficult to prove that under the same assumptions on $\varphi$
and $K$ the two conditions are equivalent:
\begin{enumerate}
\item each separable process $X(t)$, $t\in K$ which
satisfies (\ref{war2}) is a.s. bounded.
\item there exists a universal constant $S<\infty$ such that for each process $X$ the condition
(\ref{war2}) implies
\be\label{talib}
\E\sup_{s,t\in K}|X(s)-X(t)|\ls S.
\ee
\end{enumerate}
The minimal constant $S$ is denoted by $S(K,d,\varphi)$.
For a proof of this statement we refer to M. Talagrand \cite{Tal0}, Theorem 2.3.
\begin{rema}
In terms of absolutely summing operators each of the above
implications is equivalent to the fact that the injection operator
$J:\mathrm{Lip}(K)\ra C(K)$
is $(\varphi,1)$ absolutely summing, in the sense of P. Assouad, see \cite{Ass}.
\end{rema}
By far the strongest criteria for finiteness of $S(K,d,\varphi)$ were
obtained using the concept of majorizing measures which was
introduced by X. Frenique in early 70. It served him and M. Talagrand to
characterize bounded Gaussian processes. To explain briefly the concept we
introduce the following definitions.
\smallskip

\noindent
For $t\in K$ and $\va\gs 0$, we denote by  $B(t,\va)$, $S(t,\va)$ 
respectively the closed ball and the sphere
with the center at $x$ and the radius $\va$ with respect to the metric $d$, i.e. 
$$
B(t,\va):=\{s\in K:\;d(s,t)\ls\va\},\;\;S(t,\va):=\{s\in K:\;d(s,t)=\va\}.
$$
We say that $m\in \ccP(K)$ is a majorizing measure (with respect to $\varphi$ and $d$) if
$$
\ccM(m,\varphi):=\sup_{t\in K}\int^{D(K)}_0\varphi^{-1}(\frac{1}{m(B(t,\va))})d\va<\infty.
$$
X. Fernique \cite{Fer1}, \cite{Fer2} proved  that if $\varphi$ has the exponential growth then the existence
of a majorizing measure is the necessary and sufficient condition for the quantity
$S(K,d,\varphi)$ to be finite. Generalizing results of Fernique, Talagrand and
others, the author \cite{Bed} succeeded in proving that for each Orlicz function $\varphi$
the existence of a majorizing measure is always the sufficient condition for
$S(K,d,\varphi)< \infty$.
However, as it will be seen in the next chapters, the existence of a majorizing measure is not
always necessary for finiteness of $S(K,d,\varphi)$. So it is still the open problem
to characterize $(K,d)$ and $\varphi$ for which
all processes satisfying (\ref{war2}) are a.s. sample bounded.
\smallskip

\noindent
This problem was studied in depth by M. Talagrand, \cite{Tal0}.
He managed to find such a characterization of $\varphi$, $(K,d)$ in two particular,
but important for applications, cases. Namely when  $d$ is the Euclidean
distance on $\R^n$ and $K$ is a ball in $\R^n$ and the other case when $K=[-1,1]$ and
the distance $d$ is given by $d(x,y)=\eta(|x-y|)$ where $\eta$
is a concave, strictly increasing function with $\eta(0)=0$.
Generalizing his ideas and using Theorem \ref{natul1} we find the characterization in
the case when $K=T=B_{\|\cdot\|}(0,r)$ and $d(x,y)=\eta(\|x-y|\|)$
($\eta$ is concave, strictly increasing, with $\eta(0)=0$).
\smallskip

\noindent
By the definition $B(t,\va)=B_{\|\cdot\|}(t,\eta^{-1}(\va))\cap T$.
Let $\lambda$ be a normalized Lebesgue measure on $T$, that is
$\lambda(A)=\frac{|A|}{|T|}$, for each $A\in \ccB(T)$, where $|\cdot|$ is the standard Lebesgue measure on $\R^n$.
Note that
$\lambda(B(t,\va))\ls\frac{\eta^{-1}(\va)^n}{r^n}$ and
$\lambda(B(t,\va))=\frac{\eta^{-1}(\va)^n}{r^n}$ if $B_{\|\cdot\|}(t,\eta^{-1}(\va))\subset T$.
\smallskip

\noindent
The function $\eta(y)/y$ is positive and decreasing.
We assume that $\eta'(0)=\infty$ (the case of finite derivative will be considered later).
Following M. Talagrand \cite{Tal0} (Theorem 5.2) we introduce a sequence $(r_k)_{k\gs 0}$.
Let $r_{0}=\eta(r)$, for $k\gs 0$ we define
$$
r_{k+1}:=\inf\{\va\gs 0:\; r_k\ls 2\va \;\;\mbox{or}\;\;
\frac{\va}{\eta^{-1}(\va)}\ls 2\frac{r_{k}}{\eta^{-1}(r_{k})}
\}.
$$
The sequence $(r_k)_{k\gs 0}$ decreases to $0$, since $r_{k+1}\ls \frac{r_k}{2}$.
The assumption $\eta'(0)=\infty$ guarantees that $r_k>0$. There are two possibilities
$$
r_k=2 r_{k+1}
\;\;\mbox{or}\;\;
\frac{r_{k+1}}{\eta^{-1}(r_{k+1})}= 2\frac{r_{k}}{\eta^{-1}(r_{k})}.
$$
Denote by $I$ the set of $k\gs 0$ for which the first possibility holds,
and the rest by $J$. Let us notice that necessarily
\be\label{k3}
2 r_{k+1}\ls r_k,\;\; 2\frac{r_{k}}{\eta^{-1}(r_k)}\ls \frac{r_{k+1}}{\eta^{-1}(r_{k+1})}.
\ee
For $k\gs 0$ we define $\ccS_k$ as a number which satisfies the equation
$$
\int^{r_k}_{r_{k+1}}\frac{\lambda(B(0,\va))}{\va}\psi(\frac{\va}{\ccS_k\lambda(B(0,\va))})d\va=
\int^{r_k}_{r_{k+1}}\frac{\eta^{-1}(\va)^n}{r^n\va}\psi(\frac{r^n\va}{\ccS_k\eta^{-1}(\va)^n})d\va=1.
$$
If $\eta'(0)<\infty$, then there exists $m\gs 0$ such that
$r_m>0$ and $r_{m+1}=0$. That means
$$
\frac{r_m}{\eta^{-1}(r_m)}\ls \frac{\va}{\eta^{-1}(\va)}\ls 2\frac{r_m}{\eta^{-1}(r_m)},\;\;\for\;0<\va\ls r_m.
$$
We define $\ccS_m$ as the infimum over $c>0$ such that 
$$
\int^{r_m}_{0}\frac{\lambda(B(0,\va))}{\va}\psi(\frac{\va}{c \lambda(B(0,\va))})d\va=
\int^{r_m}_{0}\frac{\eta^{-1}(\va)^n}{r^n\va}\psi(\frac{r^n\va}{c \eta^{-1}(\va)^n})d\va\ls 1.
$$
For simplicity we define also $\ccS_{k}:=0$, for $k>m$ (in this case).
\smallskip

\noindent
Let us state the main result of the paper.
\bt\label{teo1}
The following inequality holds
\be\label{kwiat1}
K^{-1}\sum_{k\gs 0}\ccS_k\ls S(T,d,\varphi)\ls K\sum_{k\gs 0} \ccS_k,
\ee
where the constant $K$ depends only on $n$.
\et
In fact we show that $\sum_{k\gs 0}\ccS_k\ls 3(n+2)S(T,d,\varphi)$ and $S(T,d,\varphi)\ls (a+bn^2)\sum_{k\gs 0}\ccS_k$
for $k\gs 0$, where $a,b$ are universal constants.
\begin{coro}
Let $\varphi(x)\equiv x^p/p$, $p>1$ and $\eta(x)=x^{\alpha}$, $0<\alpha\ls 1$. Then $S(T,d,\varphi)<\infty$ 
if and only if $n<p\alpha$.
\end{coro}
In the case of $T=B_{\|\cdot\|}(0,1)$ and $d(s,t)=\|s-t\|$ (that is $\eta(x)\equiv x$)
we have $m=0$ in the above construction. Thus
$S(T,d,\varphi)$ is comparable with $\ccS_0$, where $\ccS_0$ is such that
$$
\int^{1}_{0} \va^{n-1}\psi(\frac{1}{\ccS_0\va^{n-1}})d\va=1,
$$
up to a constant which depends only on $n$. When $d$ is the Euclidean distance this corollary was proved by
M. Talagrand (Theorem 5.1, \cite{Tal0}).
\smallskip

\noindent
In the case of $T=[-1,1]$, $\eta'(0)=\infty$ it can be observed that for some universal $C>0$
$$
C^{-1}r_k\varphi^{-1}(\frac{1}{\eta^{-1}(r_k)}) \ls \ccS_k\ls Cr_k\varphi^{-1}(\frac{1}{\eta^{-1}(r_k)}),\;\;
\for\; k\gs 0.
$$
Consequently $S(T,d,\varphi)$ is comparable with $\sum_{k\gs0}r_k\varphi^{-1}(\frac{1}{\eta^{-1}(r_k)})$.
The result was obtained by M. Talagrand (Theorem 5.2, \cite{Tal0}).

\section{Preliminary results}

We remind that $\varphi,\psi$ are conjugate Orlicz functions.
\bl\label{bryndza}
Following inequalities hold:
\ba\label{nie4}
&&\varphi(\frac{\psi(x)}{x})\ls \psi(x)\ls\varphi(\frac{2\psi(x)}{x}),\;\;\for\; x\gs 0;\nonumber\\
&&\psi(\frac{\varphi(x)}{x})\ls \varphi(x)\ls\psi(\frac{2\varphi(x)}{x}),\;\;\for\; x\gs 0.
\ea
In the symmetric version we can write
$$
x\ls \varphi^{-1}(x)\psi^{-1}(x)\ls 2x,\;\;\for\;x\gs 0.
$$
\el
\begin{dwd}
Fix $x\gs 0$. By the Young inequality we obtain
$$
2\psi(x)=\frac{2\psi(x)}{x}x\ls \psi(x)+\varphi(\frac{2\psi(x)}{x}).
$$
Hence $\psi(x)\ls \varphi(\frac{2\psi(x)}{x})$. To prove the right-hand side of (\ref{nie4})
let us notice that since $\varphi(x)=\sup_{y\gs 0}(xy-\psi(y))$, we have for some $y\gs 0$
$$
\varphi(\frac{\psi(x)}{x})=y\frac{\psi(x)}{x}-\psi(y).
$$
It remains to prove that
\be\label{huta}
y\frac{\psi(x)}{x}-\psi(y)\ls\psi(x).
\ee
If $x\ls y$, then the convexity of $\psi$ gives $y\frac{\psi(x)}{x}\ls \psi(y)$.
If $x>y$, then $y\frac{\psi(x)}{x}\ls \psi(x)$. It yields (\ref{huta}),
consequently
$$
\varphi(\frac{\psi(x)}{x})\ls \psi(x)\ls\varphi(\frac{2\psi(x)}{x}).
$$
It completes the proof.
\end{dwd}
\bl\label{prr}
Functions $\varphi,\psi$ have following properties:
\begin{enumerate}
\item functions $x\varphi(1/x)$, $x\psi(1/x)$  are convex, decreasing;
\item functions $x\varphi^{-1}(1/x)$, $x\psi^{-1}(1/x)$ are concave, increasing.
\end{enumerate}
\el
\begin{dwd} It is enough to prove the result for $\varphi$.
By the definition $\varphi(x)=\sup_{y\gs 0}(yx-\psi(y))$, so $x\varphi(1/x)=\sup_{y\gs 0}(y-x\psi(y))$.
The supremum of convex functions is a convex function, the supremum of decreasing functions
is a decreasing function.
\smallskip

\noindent Similarly we observe that
$\varphi^{-1}(x)=\inf_{y\gs 0}\frac{x+\psi(y)}{y}$.
Hence $x\varphi^{-1}(1/x)=\inf_{y\gs 0}\frac{1+x\psi(y)}{y}$.
The infimum of concave functions is a concave function, the infimum of increasing functions
is an increasing function.
\end{dwd}

\section{Proof of Theorem \ref{natul1}}

\begin{dww}
Fix points $t=(t_i)^n_{i=1},\;s=(s_i)^n_{i=1} \in T$.
Let $g$ be a smooth function on $\R^n$. We define
$F_t:T\times [0,r]\ra T$ by the formula
$$
F_{t}(u,\va)=(1-\frac{\va}{r})t+\frac{\va}{r}u.
$$
We have
\be\label{nierysia}
r^n\int_T g(u)du=\int^r_0 \frac{\partial}{\partial\va}(\int_T g(F_{t}(u,\va))\va^n du)d\va.
\ee
It can be easily verified that
$$
\frac{\partial}{\partial\va}g(F_{t}(u,\va))=
r^{-1}\sum^n_{i=1}(u_i-t_i)\frac{\partial}{\partial x_i}g(F_t(u,\va))
=\va^{-1}\sum^n_{i=1}(u_i-t_i)\frac{\partial}{\partial u_i}g(F_{t}(u,\va)).
$$
Hence the following equation holds
$$
\frac{\partial}{\partial\va}(g(F_t(u,\va))\va^n)=
n\va^{n-1}g(F_{t}(u,\va))+\va^{n-1}\sum^n_{i=1}(u_i-t_i)\frac{\partial}{\partial u_i}g(F_t(u,\va)),
$$
which yields
$$
\frac{\partial}{\partial\va}(g(F_{t}(u,\va))\va^n)=
\va^{n-1}\sum^n_{i=1}\frac{\partial}{\partial u_i}(g(F_{t}(u,\va))(u_i-t_i)).
$$
Applying the generalized Green-Gauss theorem (see Theorem 4.5.6 in \cite{Fed}) which
holds for Lipschitz boundaries, we get
$$
\int_T \sum^n_{i=1}\frac{\partial}{\partial u_i}(g(F_{t}(u,\va))(u_i-t_i))du=
\int_{\partial T}g(F_{t}(u,\va))\langle u-t, n(u) \rangle \sigma_{\partial T}(du),
$$
where $\sigma_{\partial T}$ is the Lebesgue measure on the manifold $\partial T$, and
$n(u)$ the normal vector to the boundary in $u\in\partial T$, such that $\langle n(u),n(u)\rangle=1$
($n(u)$ is well defined $\sigma_{\partial T}$-a.s.).
Let us notice that the convexity of $T$ yields
$\langle u-t, n(u) \rangle\gs 0$.
Denoting $\sigma_t(u):=\langle u-t, n(u) \rangle\sigma_{\partial T}(du)$, we obtain due to (\ref{nierysia}) 
$$
r^n\int_T g(u)du=
\int^r_0\int_{\partial T}g(F_t(u,\va)) \va^{n-1}\sigma_t(du)d\va.
$$
By the standard approximation this equality can be easily generalized to any Borel, bounded
function $g$ on $T$. We verify also that $n|T|=\sigma_t(\partial T)$ (consider $g\equiv 1$), 
consequently for each $g\in \mathfrak{B}(T)$
\be\label{nie1}
n^{-1}r^n\pint_T g(u)du=
\int^r_0\pint_{\partial T}g(F_t(u,\va)) \va^{n-1}\sigma_t(du)d\va. 
\ee
We define $\partial F_t:\partial T\times [0,1]\ra T$ by 
$\partial F_t=F_{t}| \partial T\times [0,1]$.
The equation (\ref{nie1}) implies 
\be\label{pis1}
|A|/|T|=\sigma_t\otimes d(\va/r)^n((\partial F_t)^{-1}(A)),\;\; \for\; A\in \ccB(T).
\ee
\smallskip

\noindent
Let $a_t:[0,r]\ra \R$ denotes
$$
a_t(\va):=\pint_{\partial T}f(F_t(u,\va))\sigma_t(du)=
\pint_{\partial T}f((1-\frac{\va}{r})t+\frac{\va}{r}u)\sigma_t(du).
$$
Clearly $a_t$ satisfies Lipschitz condition (because $f$ does) and
$$
a_t(0)=f(t),\;\;a_t(r)=\pint_{\partial T} f(u)\sigma_t(du).
$$
Since $f$ is Lipschitz, there exists bounded $\nabla f$, $|\cdot|$-a.s. on $T$.
We check that if $f$ is differentiable in $F_t(u,\va)$ then
$$
\frac{\partial}{\partial \va}f(F_t(u,\va))=r^{-1}\langle\nabla f(F_t(u,\va)), u-t \rangle.
$$
By the Fubini theorem and (\ref{pis1}) we have that $f(F_t)$ is differentiable $d\va$-a.s. for 
$\sigma_t(du)$-almost all $u\in \partial T$.
Consequently $d\va$-a.s. there holds
$$
a_t'(\va)=r^{-1}\pint_{\partial T}\langle \nabla f(F_t(u,\va)), u-t \rangle \sigma_t(du).
$$
We have
\ba
&&|a_t'(\va)|=r^{-1}|\pint_{\partial T}\langle \nabla f(F_t(u,\va)), u-t \rangle \sigma_t(du)|\ls\nonumber\\
&&\ls r^{-1}\pint_{\partial T}\|u-t\|
\|\nabla f(F_t(u,\va))\|_{\ast}\sigma_t(du).\nonumber
\ea
Clearly $\|u-t\|\ls 2r$. Observe that $b_t:[0,r]\ra \R$ is $d\va$-a.s. well defined by
the formula
$$
b_t(\va):=\pint_{\partial T}\|\nabla f(F_t(u,\va))\|_{\ast}\sigma_t(du).
$$
Hence $|a_t'(\va)|\ls 2b_t(\va)$, $d\va$-a.s. By the Jensen inequality
and (\ref{nie1}) for each $B>0$ we obtain 
\ba
&&\int^{r}_0 \varphi(\frac{1}{B}b_t(\va))\va^{n-1}d\va\ls \pint_{\partial T}\int^{r}_0
\varphi(\frac{1}{B}
\|\nabla f(F_t(u,\va))\|_{\ast})
\va^{n-1}d\va\sigma_t(du)
=\nonumber\\
&&=n^{-1}r^{n}\pint_T\varphi(\frac{1}{B}\|\nabla f(u)\|_{\ast})du.\nonumber
\ea
The Young inequality (\ref{nie3})
gives
$$
\frac{b_t(\va)}{AB\va^{n-1}}\ls \psi(\frac{1}{A\va^{n-1}})+\varphi(\frac{1}{B}b_t(\va)).
$$
Since $a_t$ is Lipschitz we get
$$
|f(t)-\pint_{\partial T}f(u)\sigma_t(u)|=|a_t(0)-a_t(r)|=|\int^{r}_0 a_t'(\va)d\va|\ls 2\int^r_0 b_t(\va)d\va.
$$
Thus
\ba\label{nie2}
&&|f(t)-\pint_{\partial T}f(u)\sigma_t(u)|\ls\nonumber\\
&&\ls 2AB(
\int^r_0 \psi(\frac{1}{A\va^{n-1}})\va^{n-1}d\va+
n^{-1}r^n\pint_T\varphi(\frac{1}{B}\|\nabla f(u)\|_{\ast})du).
\ea
Again due to the generalized Green-Gauss theorem (this version
holds for Lipschitz functions and Lipschitz boundaries) we obtain 
\ba
&&|\pint_{\partial T}f(u)\sigma_t(du)-\pint_{\partial T}f(u)\sigma_s(du)|=\nonumber\\
&&=\frac{1}{n|T|}|\int_{\partial T}f(u)\langle s-t,n(u) \rangle \sigma_{\partial T}(du)|=\nonumber\\
&&=n^{-1}|\pint_T \langle \nabla f(u),s-t\rangle du|
\ls 2rn^{-1} \pint_T \|\nabla f(u)\|_{\ast}du.\nonumber
\ea
By the Young inequality and since $y\psi(1/y)$ is decreasing we obtain
\ba
&&\frac{\|\nabla f(u)\|_{\ast}}{ABr^{n-1}}\ls
\psi(\frac{1}{Ar^{n-1}})+\varphi(\frac{1}{B}\|\nabla f(u)\|_{\ast})\ls\nonumber\\
&&\ls n\int^r_0 \psi(\frac{1}{A\va^{n-1}})\frac{\va^{n-1}}{r^n}d\va+\varphi(\frac{1}{B}\|\nabla f(u)\|_{\ast}).\nonumber
\ea
It follows that
\ba\label{jp2}
&&|\pint_{\partial T}f(u)\sigma_t(du)-\pint_{\partial T}f(u)\sigma_s(du)|\ls\nonumber\\
&&\ls 2AB(\int^r_0 \psi(\frac{1}{A\va^{n-1}})\va^{n-1}d\va+
n^{-1}r^n\pint_T\varphi(\frac{1}{B}\|\nabla f(u)\|_{\ast})du).
\ea
By the definition $|T|=|B(0,r)|=r^n|B_{\|\cdot\|}(0,1)|$.
Inequalities (\ref{nie2}), (\ref{jp2}) yield
\ba
&&|f(s)-f(t)|\ls |f(t)-\pint_{\partial T}f(u)\sigma_t(u)|+\nonumber\\
&&+|\pint_{\partial T}f(u)\sigma_t(du)-\pint_{\partial T}f(u)\sigma_s(du)|+|f(s)-\pint_{\partial T}f(u)\sigma_s(u)|\ls \nonumber\\
&&\ls 6AB(\int^r_0 \psi(\frac{1}{A\va^{n-1}})\va^{n-1}d\va+
\frac{1}{n|B_{\|\cdot\|}(0,1)|}\int_T\varphi(\frac{1}{B}\|\nabla f(u)\|_{\ast} )du).\nonumber
\ea
\end{dww}
\begin{rema}\label{truskawka}
Let $\sigma_{\partial T}$ be the Lebesgue measure on the manifold $\partial T$ and
$n(u)$ the normal vector to the boundary in $u\in \partial T$ such that $\langle n(u),n(u)\rangle=1$ ($n(u)$
is well defined $\sigma_{\partial T}$-a.s.).
For each $t\in T$ and measure $\sigma_t(du)=\langle u-t, n(u) \rangle \sigma_{\partial T}(du)$
the equality $\sigma_{t}(\partial T)=n|T|=nr^n|B_{\|\cdot\|}(0,1)|$ holds.
\end{rema}
\begin{coro}\label{gry0}
For each $t\in T$, there holds
\ba
&&|f(t)-\pint_{T}f(u)du|\ls 6AB
(\int^r_{0} \psi(\frac{1}{A\va^{n-1}})\va^{n-1}d\va
+\nonumber\\
&&+\frac{1}{n|B_{\|\cdot\|}(0,1)|}\int_T\varphi(\frac{1}{B}\|\nabla f(u)\|_{\ast} )du).\nonumber
\ea
\end{coro}
\begin{dwd}
This observation is obvious.
Theorem \ref{natul1} yields
\ba
&&|f(t)-f(u)|\ls 6AB
(\int^r_0 \psi(\frac{1}{A\va^{n-1}})\va^{n-1}d\va+\nonumber\\
&&+\frac{1}{n|B_{\|\cdot\|}(0,1)|}\int_T\varphi(\frac{1}{B}\|\nabla f(u)\|_{\ast} )du).\nonumber
\ea
Integrating both parts and using
$|f(t)-\pint_{T}f(u)du|\ls\pint_T |f(t)-f(u)|du$ we obtain the corollary.
\end{dwd}

\section{Construction of the optimal process}

We assume that $\eta'(0)=\infty$.
In this section we prove the left-hand side of (\ref{kwiat1}) in Theorem \ref{teo1}.
\begin{dw1}
We define a stochastic process on a probability space
$(T,\ccB(T),\lambda)$ by the formula
$$
X(t,\omega):=\int^{d(\omega,t)}_{d(\omega,0)}g(\va)d\va,\;\;\for
\;t\in T,\omega\in T,
$$
where $g(\va)$ is a positive function, integrable on each interval
$[\delta,\eta(r)]$, $\delta>0$ and such that $g(\va)=0$, for $\va> \eta(r)$.
Let us notice that the process $X$ is separable. Suppose we have shown that
$$
\E\varphi(\frac{|X(s)-X(t)|}{d(s,t)})=
\int_T \varphi(\frac{|X(s,\omega)-X(t,\omega)|}{d(s,t)})\lambda(d\omega)\ls 1,
$$
then the process $X(t)$, $t\in T$ satisfies the condition (\ref{war2}). Since
$$
X(\omega,\omega)=-\int^{d(\omega,0)}_{0}g(\va)d\va,\;\;
X(\frac{\|\omega\|-r}{\|\omega\|}\omega,\omega)=\int^{\eta(r)}_{d(\omega,0)} g(\va)d\va,
$$
we have $\sup_{s,t\in T}|X(s,\omega)-X(t,\omega)|=\int^{\eta(r)}_0 g(\va)d\va$.
Due to the definition of $S(T,d,\varphi)$ it proves that
\be\label{nat1}
\int^{\eta(r)}_0 g(\va)d(\va)\ls \E\sup_{s,t\in T}|X(s)-X(t)|\ls S(T,d,\varphi).
\ee
The convexity of $\varphi$, $\varphi(0)=0$
and the Jensen inequality imply
\ba
&&\int_T \varphi(|\frac{X(s,\omega)-X(t,\omega)|}{d(s,t)})\lambda(d\omega)\ls\nonumber\\
&&\ls \int_T\varphi(\frac{|d(s,\omega)-d(t,\omega)|}{d(s,t)}
\pint^{d(t,\omega)}_{d(s,\omega)}g(\va)
d\va ) \lambda(d\omega)\ls \nonumber\\
&&\ls
\int_T\frac{|d(s,\omega)-d(t,\omega)|}{d(s,t)}
|\pint^{d(t,\omega)}_{d(s,\omega)}\varphi(g(\va))
d\va| \lambda(d\omega)=\nonumber\\
&&=\frac{1}{d(s,t)}\int_T|\int^{d(t,\omega)}_{d(s,\omega)}\varphi(g(\va))
d\va| \lambda(d\omega).\nonumber
\ea
The Fubini theorem yields
$$
\int_T|\int^{d(t,\omega)}_{d(s,\omega)}\varphi(g(\va))
d\va| \lambda(d\omega)=
\int^{\eta(r)}_0\varphi(g(\va))\lambda(B(s,\va)\triangle B(t,\va))d\va,
$$
where $\triangle$ is the symmetric set difference.
Observe that if $d(s,t)\gs \va$, then
$$
\lambda(B(s,\va)\triangle B(t,\va))\ls
\lambda(B(s,\va))+\lambda(B(t,\va))\ls 2\frac{\eta^{-1}(\va)^n}{r^n}.
$$
From the other hand if $\va\gs d(s,t)$, then
$$
B_{\|\cdot\|}(\frac{s+t}{2},\eta^{-1}(\va)-\frac{1}{2}\|s-t\|)\subset
B_{\|\cdot\|}(s,\eta^{-1}(\va))\cap B_{\|\cdot\|}(t,\eta^{-1}(\va)),
$$
and thus
\ba
&&\frac{|B(s,\va)\triangle B(t,\va)|}{|B(0,\eta(r))|}\ls 2(\frac{\eta^{-1}(\va)^n}{r^n}-\frac{(\eta^{-1}(\va)-\frac{1}{2}\|s-t\|)^n}{r^n})\ls\nonumber\\
&&\ls n\|s-t\|\frac{\eta^{-1}(\va)^{n-1}}{r^n}.\nonumber
\ea
Hence, for $\va\gs d(s,t)$ we have $\lambda(B(s,\va)\triangle B(t,\va))\ls n\|s-t\|\frac{\eta^{-1}(\va)^{n-1}}{r^n}$.
Consequently if $d(s,t)\gs \eta(r)$ then
\ba\label{harpie0}
&&\int_T \varphi(\frac{|X(s,\omega)-X(t,\omega)|}{d(s,t)})\mu(d\omega)\ls
\frac{2}{d(s,t)}\int^{\eta(r)}_0\frac{\eta^{-1}(\va)^n}{r^n}\varphi(g(\va))d\va
\ls\nonumber\\
&&\ls \frac{2}{\eta(r)}\int^{\eta(r)}_0\frac{\eta^{-1}(\va)^n}{r^n}\varphi(g(\va))d\va
\ea
and if $d(s,t)\ls \eta(r)$, then
\ba\label{harpie}
&&\int_T \varphi(|\frac{X(s,\omega)-X(t,\omega)|}{d(s,t)})\lambda(d\omega)\ls
\frac{2}{d(s,t)}\int^{d(s,t)}_0\frac{\eta^{-1}(\va)^n}{r^n}\varphi(g(\va))d\va+\nonumber\\
&&+\frac{\|s-t\|}{d(s,t)}\int^{\eta(r)}_{d(s,t)}n\frac{\eta^{-1}(\va)^{n-1}}{r^n}\varphi(g(\va))d\va.
\ea
The construction of $g$ is as follows
$$
g(\va):=K^{-1}\frac{\ccS_k\eta^{-1}(\va)^{n}}{r^n\va}\psi(\frac{r^n\va}{\ccS_k\eta^{-1}(\va)^n}),\;\;\for\; r_{k+1}<\va\ls r_k,
$$
where the constant $K\gs 1$ we choose later. From the convexity of $\varphi$ and Lemma \ref{bryndza} we deduce
\be\label{k2}
\varphi(g(\va))\ls K^{-1}\psi(\frac{r^n\va}{\ccS_k\eta^{-1}(\va)^n}),\;\;\for\;r_{k+1}<\va\ls r_k,\; k\gs 0.
\ee
We show that the process $X$ satisfies the condition (\ref{war2}) for such $g$.
\smallskip

\noindent
First we assume that $d(s,t)\ls \eta(r)=r_0$. Hence there exists $m$ such that $r_{m+1}<d(s,t)\ls r_m$. Consider $k>m$,
the definition of $\ccS_k$ and (\ref{k2}) yield
$$
\int^{r_k}_{r_{k+1}}\frac{\eta^{-1}(\va)^n}{r^n}\varphi(g(\va))d\va\ls
K^{-1}r_{k}\int^{r_k}_{r_{k+1}}\frac{\eta^{-1}(\va)^n}{r^n\va}\psi(\frac{r^n\va}{\ccS_k\eta^{-1}(\va)^n})d\va\ls K^{-1}r_k.
$$
Similarly we obtain
$$
\int^{d(s,t)}_{r_{m+1}}\frac{\eta^{-1}(\va)^n}{r^n}\varphi(g(\va))d\va\ls
K^{-1}d(s,t).
$$
Since (\ref{k3}) gives $2r_{k+1}\ls r_k$, it is clear that $r_k\ls 2^{-k+m+1}d(s,t)$, for $k>m$.
Applying the above inequalities, we get
\be\label{nyt}
\frac{2}{d(s,t)}\int^{d(s,t)}_0\frac{\eta^{-1}(\va)^n}{r^n}\varphi(g(\va))d\va
\ls 2K^{-1}(1+\sum_{k>m}2^{-k+{m+1}})= 6K^{-1}.
\ee
It remains to find the estimation for the second integral in (\ref{harpie}).
Consider $0\ls k< m$. By (\ref{k2}), the definition of $\ccS_k$ and since $\eta^{-1}(y)/y$ is 
increasing, we have
\ba\label{tr2}
&&\int^{r_k}_{r_{k+1}}n\frac{\eta^{-1}(\va)^{n-1}}{r^n}\varphi(g(\va))d\va\ls
nK^{-1}\int^{r_k}_{r_{k+1}}\frac{\eta^{-1}(\va)^{n-1}}{r^n}\psi(\frac{r^n\va}{\ccS_k\eta^{-1}(\va)^n})d\va\ls
\nonumber\\
&&
\ls nK^{-1}\frac{r_{k+1}}{\eta^{-1}(r_{k+1})}
\int^{r_k}_{r_{k+1}}\frac{\eta^{-1}(\va)^n}{r^n\va}\psi(\frac{r^n\va}{\ccS_k\eta^{-1}(\va)^n})d\va=
\nonumber\\
&&= nK^{-1}\frac{r_{k+1}}{\eta^{-1}(r_{k+1})}.
\ea
In the same way, we prove
\be\label{tr3}
\int^{r_m}_{d(s,t)}n\frac{\eta^{-1}(\va)^{n-1}}{r^n}\varphi(g(\va))d\va\ls nK^{-1}\frac{d(s,t)}{\|s-t\|}.
\ee
Let us notice that (\ref{k3}) gives
$$
\frac{\|s-t\|}{d(s,t)}\frac{r_{k+1}}{\eta^{-1}(r_{k+1})}\ls 2^{-m+k+1},\;\;\for\;0\ls k<m.
$$
Inequalities (\ref{tr2}) and (\ref{tr3}) imply
\ba\label{nat}
&&\frac{\|s-t\|}{d(s,t)}\int^{\eta(r)}_{d(s,t)}n\frac{\eta^{-1}(\va)^{n-1}}{r^n}\varphi(g(\va))d\va\ls\nonumber\\
&&\ls nK^{-1}(1+\sum^{m-1}_{k=0}2^{-m+k+1})\ls 3nK^{-1}.
\ea
If we plug estimations (\ref{nyt}), (\ref{nat}) into (\ref{harpie}), we obtain
$$
\E\varphi(\frac{|X(s)-X(t)|}{d(s,t)})\ls 3K^{-1}(2+n).
$$
The second case is when $d(s,t)\gs \eta(r)$. We use (\ref{nyt}) to get
$$
\frac{2}{\eta(r)}\int^{\eta(r)}_0\frac{\eta^{-1}(\va)^n}{r^n}\varphi(g(\va))d\va
\ls 6K^{-1}.
$$
The above inequality and (\ref{harpie0}) imply 
$$
\E\varphi(\frac{|X(s)-X(t)|}{d(s,t)})\ls 6K^{-1}\ls 3K^{-1}(2+n).
$$ 
Therefore for (\ref{war2}) we need that $K:=3(2+n)$.
\smallskip

\noindent
By the definition of $g$ and numbers $\ccS_k$
it is clear that for all $k\gs 0$ we have $\int^{r_k}_{r_{k+1}}g(\va)d\va=K^{-1}\ccS_k$.
Consequently (\ref{nat1}) gives
$\sum^{\infty}_{k=0}\ccS_k\ls KS(T,d,\varphi)$. The theorem is proved with the constant $K=3(2+n)$.
\end{dw1}

\section{Some basic tools}

Before we prove the right-hand side of (\ref{kwiat1}) we establish
some helpful results. We start from proving a fact which allows us to
consider processes with finite number of different Lipschitz paths.
\begin{lema}\label{jujki}
Fix any point $t_0\in T$. Let $F\subset T$ be a finite set.
For each process $X(t)$, $t\in T$ which satisfies the condition (\ref{war2}) there exists
a sequence of processes $(Y_k)_{k\gs 1}$ which  satisfy (\ref{war2}), have
finite number of different Lipschitz trajectories and such that
\be\label{bre}
\lim_{k\ra \infty}Y_k(t)=X(t)-X(t_0),\;\;\mbox{a.s. and  in}\; L_{1},\;\; \for\;t\in F.
\ee
In particular (\ref{bre}) implies
$$
\lim_{k\ra\infty}\E\sup_{s,t\in F}|Y_k(s)-Y_k(t)|=\E\sup_{s,t\in F}|X(s)-X(t)|.
$$
\end{lema}
\begin{dwd}
A process $X(t)$, $t\in T$ is defined on a probability space $(\Omega,\ccF,\P)$.
Denote $Y(t):=X(t)-X(t_0)$. It is clear that $Y(t)$, $t\in T$
satisfies condition (\ref{war2}) and moreover $\|Y(t)\|_{\varphi}\ls d(t,t_0)$,
what implies $\E|Y(t)|<\infty$, for $t\in T$.
\smallskip

\noindent
First we assume that $\ccF$ is a finite $\sigma$-algebra.
Due to (\ref{war2}) we have
$$
|Y(s,\omega)-Y(t,\omega)|\ls d(s,t)\varphi^{-1}(1/\P(A)),\;\;\for\;s,t\in T,\;\omega\in A,
$$
where $A$ is an atom in $\ccF$.
Hence the process $Y$ has $\P$-a.s. finite number of different Lipschitz trajectories.
\smallskip

\noindent
In the general case we use the fact that $F$ is a finite set.
There exists an increasing sequence of finite $\sigma$-algebras $(\ccF_k)_{k\gs 1}$
which sum generates $\sigma(Y(t):\;t\in F)$.  Notice that $\E|Y(t)|<\infty$, for $t\in T$ 
(since $\|Y(t)\|_{\varphi}<\infty$), thus we can define $Y_k(t):=\E(Y(t)|\ccF_k)$, $t\in T$. By the Jensen inequality we get
$$
\E\varphi(\frac{|Y_k(s)-Y_k(t)|}{d(s,t)})\ls\E\varphi(\frac{|Y(s)-Y(t)|}{d(s,t)})\ls 1,\;\;\for\;s,t\in T.
$$
The process $Y_k$ satisfies (\ref{war2}), hence $\P$-a.s. it has finite number of different 
Lipschitz trajectories. Modifying $Y_k$ on the set of measure $0$ we may assume that $Y_k$
has finite number of different Lipschitz trajectories.
Clearly $Y_k(t)\ra Y(t)$ $\P$-a.s., for $t\in F$. Since $\E|Y(t)|<\infty$, the convergence is also in
$L_{1}$.
\end{dwd}
Next step is to prove some approximation on numbers $\ccS_k$.
\bl\label{bromba}
There holds:
\ba
&&\frac{1}{4}r_{k}\varphi^{-1}(\frac{2r^n}{\eta^{-1}(r_{k})^n})\ls \ccS_k,\;\;\for\; k\gs 0;\nonumber\\
&&\ccS_k\ls r_{k+1}\varphi^{-1}(\frac{r^n}{\eta^{-1}(r_{k+1})^n}),\;\;\for\; k\in I.\nonumber
\ea
\el
\begin{dwd}
Due to (\ref{k3}) we know that $r_k-r_{k+1}\gs \frac{1}{2}r_k$.
Lemma \ref{prr} follows that $y\psi(1/y)$ is decreasing. Thus, for $k\gs 0$ we have
\ba
&&1=\int^{r_k}_{r_{k+1}}\frac{\eta^{-1}(\va)^n}{r^n\va}\psi(\frac{r^n\va}{\ccS_k\eta^{-1}(\va)^n})d\va\gs
\nonumber\\
&&\gs(r_k-r_{k+1})\frac{\eta^{-1}(r_k)^n}{r^nr_k}\psi(\frac{r^nr_k}{\ccS_k\eta^{-1}(r_k)^n})\gs
\frac{\eta^{-1}(r_{k})^n}{2r^n}\psi(\frac{r^n r_{k}}{\ccS_k\eta^{-1}(r_{k})^n}).\nonumber
\ea
That means
$$
\psi^{-1}(\frac{2r^n}{\eta^{-1}(r_{k})^n})\gs \frac{r^n r_{k}}{\ccS_k\eta^{-1}(r_{k})^n}.
$$
By Lemma \ref{bryndza} (that is by the inequality $\varphi^{-1}(y)\psi^{-1}(y)\ls 2y$) we obtain
$$
\ccS_k\gs \frac{1}{4}r_{k}\varphi^{-1}(\frac{2r^n}{\eta^{-1}(r_{k})^n}).
$$
We prove the second inequality. Since $y\psi(1/y)$ is decreasing and $r_k-r_{k+1}=r_{k+1}$, for $k\in I$, then
\ba
&&1=\int^{r_k}_{r_{k+1}}\frac{\eta^{-1}(\va)^n}{r^n \va}\psi(\frac{r^n \va}{\ccS_k\eta^{-1}(\va)^n})d\va\ls\nonumber\\
&&
\ls r_{k+1}\frac{\eta^{-1}(r_{k+1})^n}{r_{k+1}r^n}\psi(\frac{r^nr_{k+1}}{\ccS_k \eta^{-1}(r_{k+1})^n})=\frac{\eta^{-1}(r_{k+1})^n}{r^n}
\psi(\frac{r^n r_{k+1}}{\ccS_k\eta^{-1}(r_{k+1})^n}).\nonumber
\ea
Hence
$$
\psi^{-1}(\frac{r^n}{\eta^{-1}(r_{k+1})^n})\ls \frac{r^n r_{k+1}}{\ccS_k\eta^{-1}(r_{k+1})^n}.\nonumber
$$
Again, using Lemma \ref{bryndza} (the inequality $y\ls \varphi^{-1}(y)\psi^{-1}(y)$), we get
$$
\ccS_k\ls r_{k+1}\varphi^{-1}(\frac{r^n}{\eta^{-1}(r_{k+1})^n}).
$$
\end{dwd}
Let us remind  that $\lambda$ is the normalized Lebesgue measure on $T$.
For $0<\va\ls \eta(r)$, we denote 
$$
B_{\va}(t):=B((1-\frac{\eta^{-1}(\va)}{r})t,\va),\;\;
S_{\va}(t):=S((1-\frac{\eta^{-1}(\va)}{r})t,\va).
$$
Observe that $B_{\va}(t)=B_{\|\cdot\|}((1-\frac{\eta^{-1}(\va)}{r})t,\eta^{-1}(\va))\subset T$, hence $\lambda(B_{\va}(t))=\frac{\eta^{-1}(\va)^n}{r^n}$. For each $f\in C(T)$ we define $f_{\va}(t):=\pint_{B_{\va}(t)}f(u)\lambda(du)$.
\smallskip

\noindent
Let us assume that $0<\va\ls \eta(r)$.
We denote by $\sigma_{t,\va}$ the Lebesgue measure on the manifold
$S_{\va}(t)$. For each $e\in\R^n$ we define
$$
\triangle^e_{t,\va}:=\{u\in S_{\va}(t):\; \langle e, n(u)\rangle \gs 0  \},
$$
where $n(u)$ is the normal vector in $u\in S_{\va}(t)$,
such that $\langle n(u),n(u) \rangle =1$ ($n(u)$ is well defined $\sigma_{t,\va}$-a.s.).
Observe that for $e\in\R^n$ and $f\in C(T)$
\ba\label{grwat}
&&\lim_{h\ra +0}\frac{1}{h}\int_{B_{\va}(t+he)\backslash B_{\va}(t)}f(u)du=
\lim_{h\ra +0}\frac{1}{h}\int_{B_{\va}(t)\backslash B_{\va}(t-he)}f(u)du=\nonumber\\
&&=(1-\frac{\eta^{-1}(\va)}{r})\int_{\triangle^e_{t,\va}}
f(u)\langle e, n(u)\rangle\sigma_{t,\va}(du).
\ea
Let $\sigma^e_{t,\va}$ denotes the positive measure on $\triangle^e_{t,\va}$ given by the formula 
$\sigma^e_{t,\va}(du):=\langle e, n(u)\rangle\sigma_{t,\va}(du)$.
Notice that if $f\in \mathrm{Lip}(T)$, then there exists $\nabla f$, $|\cdot|$-a.s.  
\bl\label{bagi}
Fix $0<\va\ls \eta(r)$. For each $e\in \R^n$,
$\|e\|=1$ and $f\in \mathrm{Lip}(T)$ the following equality
holds $|\cdot|$-a.s. on $T$
$$
\langle \nabla f_{\va}(t),e \rangle
=(1-\frac{\eta^{-1}(\va)}{r}) \beta(e) \pint_{\triangle^{e}_{t,\va}}\pint_{\triangle^{-e}_{t,\va}}\frac{f(u)-f(v)}{2\eta^{-1}(\va)}
\sigma^{e}_{t,\va}(du)\sigma^{-e}_{t,\va}(dv),
$$
where $\beta(e)\ls n$.
\el
\begin{dwd}
First we assume $h>0$. Observe that
\ba\label{watabi}
&&f_{\va}(t+he)-f_{\va}(t)=\pint_{B_{\va}(t+he)}f(u)\lambda(du)-\pint_{B_{\va}(t)}f(u)\lambda(du)=\nonumber\\
&&=\frac{1}{|B_{\va}(t)|}
(\int_{B_{\va}(t+he)\backslash B_{\va}(t)}f(u)du-\int_{B_{\va}(t)\backslash B_{\va}(t+he)}f(u)du).
\ea
By (\ref{grwat}) we obtain
\be\label{lisek1}
\lim_{h\ra +0}\frac{1}{h}\int_{B_{\va}(t+he)\backslash B_{\va}(t)}f(u)du=(1-\frac{\eta^{-1}(\va)}{r}) \int_{\triangle^{e}_{t,\va}}f(u)
\sigma^{e}_{t,\va}(du).
\ee
Let us define $\beta(e)$ by the following formula
$$
\beta(e):=\frac{2\eta^{-1}(\va)\sigma^{e}_{t,\va}(\triangle^{e}_{t,\va})}{|B_{\va}(t)|}
=\frac{2\eta^{-1}(\va)\int_{\triangle^{e}_{t,\va}}\langle e,n(u)\rangle\sigma_{t,\va}(du) }{|B_{\va}(t)|}.
$$
The homogeneity and symmetry imply
$$
\beta(e)=
\frac{\int_{\triangle^{e}_{0,\eta(r)}}\langle r e,n(u)\rangle\sigma_{0,\eta(r)}(du)}{|B_{\eta(r)}(0)|}
+\frac{\int_{\triangle^{-e}_{0,\eta(r)}}\langle -r e,n(u)\rangle\sigma_{0,\eta(r)}(du)}{|B_{\eta(r)}(0)|}.
$$
Due to Remark \ref{truskawka}, for each $t\in B_{\eta(r)}(0)=B(0,\eta(r))$ we have
$$
\int_{S_{\eta(r)}(0)}\langle u-t,n(u)\rangle\sigma_{0,\eta(r)}(du)=n|B_{\|\cdot\|}(0,r)|=n|B_{\eta(r)}(0)|.
$$
Applying this equality for $t=-r e$, $t=re$ and $t=0$, we get
\ba
&&\beta(e)=\frac{\int_{\triangle^{e}_{0,\eta(r)}}\langle u+r e,n(u)\rangle\sigma_{0,\eta(r)}(du)}{|B_{\eta(r)}(0)|}
+\nonumber\\
&&+\frac{\int_{\triangle^{-e}_{0,\eta(r)}}\langle u-r e,n(u)\rangle\sigma_{0,\eta(r)}(du)}{|B_{\eta(r)}(0)|}-\frac{\int_{\triangle^{e}_{0,\eta(r)}}\langle u,n(u)\rangle\sigma_{0,\eta(r)}(du)}{|B_{\eta(r)}(0)|}- \nonumber\\
&&-\frac{\int_{\triangle^{-e}_{0,\eta(r)}}\langle u,n(u)\rangle\sigma_{0,\eta(r)}(du)}{|B_{\eta(r)}(0)|}
\ls (n+n-n)=n.\nonumber
\ea
We have used here the fact that
$\sigma_{0,\eta(r)}( S_{\eta(r)}(0)\backslash(\triangle^{e}_{0,\eta(r)}\cup \triangle^{-e}_{0,\eta(r)}))=0$.
Observe that by (\ref{grwat})
\be\label{lisek2}
\lim_{h\ra +0}\frac{1}{h}\int_{B_{\va}(t)\backslash B_{\va}(t+he) }f(u)\lambda(du)=
(1-\frac{\eta^{-1}(\va)}{r})\int_{\triangle^{-e}_{t,\va}}f(u)
\sigma^{-e}_{t,\va}(du).
\ee
Moreover $\beta(e)=\beta(-e)$, hence applying (\ref{watabi}), (\ref{lisek1}) and (\ref{lisek2}) we obtain
\ba
&&\lim_{h\ra +0}\frac{1}{h}(f_{\va}(t+he)-f_{\va}(t))=\nonumber\\
&&=(1-\frac{\eta^{-1}(\va)}{r})\beta(e)
\pint_{\triangle^{e}_{t,\va}}\pint_{\triangle^{-e}_{t,\va}}\frac{f(u)-f(v)}{2\eta^{-1}(\va)}
\sigma^{e}_{t,\va}(du)
\sigma^{-e}_{t,\va}(dv).\nonumber
\ea
The case of $h<0$ can be treated in the similar way.
\end{dwd}

\section{The estimation from above}

We assume that $\eta'(0)=\infty$.
In this section we prove the right-hand side of (\ref{kwiat1}) in Theorem \ref{teo1}.
\begin{dw2} Denote $B_k(x):=B_{r_k}(x)$, $\for\;x\in T$. Let us notice that
\be\label{kluski}
\lambda(B_k(x))=\lambda(B_{r_k}(x))=\frac{\eta^{-1}(r_k)^n}{r^n},\;\;\for\;x\in T.
\ee
For each $k\gs 0$ we define a linear operator
$S_k:C(T)\ra C(T)$ by the formula
$$
S_k f(x):=f_{r_k}(x)=\pint_{B_{k}(x)}f(u)\lambda(du),\;\;\for\; x\in T.
$$
If $f,g\in C(T)$, $k\gs 0$, then such properties can be easily derived:
\begin{enumerate}
\item $S_k 1=1$;
\item if $f\ls g$, then $S_k f\ls S_k g$ and so $|S_k f|\ls S_k |f|$;
\item $S_{0} f=\int_T f(u)\lambda(du)$ hence $S_{k}S_{0} f=S_{0} f$;
\item $\lim_{k\ra\infty} S_k f(x)=f(x)$.
\end{enumerate}
Observe that if $f\in \mathrm{Lip}(T)$, then
$S_k f$ is also Lipschitz and thus differentiable $|\cdot|$-a.s.
Fix $m>0$. For $0\ls k\ls m$, we define operators
$L_k,R_k,T_k:C(T)\ra C(T)$. We put $L_m=Id,R_m=Id$ and for $k<m$ 
\ba
&&L_k:=L_{k+1}S_{k+1},\;\;R_k:=R_{k+1},\;\;\mbox{if }\;k\in I;\nonumber\\
&&L_k:=L_{k+1},\;\;R_k:=S_{k+1}R_{k+1},\;\;\mbox{if}\;k\in J.\nonumber
\ea
Denote also $T_{k}:=L_kS_kR_k$, for $0\ls k\ls m$.
Properties of $S_k$ imply corresponding
properties of operators $L_k,R_k,T_k$.
For $0\ls k\ls m$ and $f,g\in C(T)$, we have:
\begin{enumerate}
\item $L_k 1=1$, and if $f\ls g$, then $L_k f\ls L_k g$;
\item if $f\in \mathrm{Lip}(T)$, then $R_k f$ is $|\cdot|$-a.s. differentiable on $T$;
\item the function $T_0 f$ is constant;
\item there holds $|T_m f(t)-T_0 f(t)|\ls \sum^{m-1}_{k=0}|T_k f(t)-T_{k+1} f(t)|$.
\end{enumerate}
Fix $f\in \mathrm{Lip}(T)$ and points $s,t\in T$.
We will analyse $|T_{k+1}f(t)-T_k f(t)|$.
There are two cases. Either $k\in I$ or $k\in J$.
In fact we use two different methods.
\smallskip

\noindent {\bf Case 1.} Fix $k\in I$, $k<m$.
By the definition, we have
$$
T_{k+1}f(t)-T_k f(t)=L_{k+1}S_{k+1}(Id-S_k)R_{k+1}f(t).
$$
Clearly $R_k=R_{k+1}$, for $k\in I$. Denote $g:=R_{k}f$,
it can be easily checked
$$
|S_{k+1}(Id-S_k)g(w)|\ls \pint_{B_{k+1}(w)}\pint_{B_k(u)}|g(u)-g(v)|\lambda(dv)\lambda(du).
$$
For each Orlicz function $\varphi$ there holds
$$
x\ls 1+\frac{\varphi(xy)}{\varphi(y)},\;\;\for\;x\gs 0,y>0.
$$
Thus
$$
\frac{|g(u)-g(v)|}{10r_k\varphi^{-1}(\frac{1}{\lambda(B_{k+1}(w))})}\ls
1+\lambda(B_{k+1}(w))
\varphi(\frac{|g(u)-g(v)|}{10r_k}).
$$
Consequently, for $u\in B_{k+1}(w)$ the inequality holds
$$
|g(u)-g(v)|\ls 10r_{k}\varphi^{-1}(\frac{1}{\lambda(B_{k+1}(w))})(1+
\lambda(B_{k+1}(w))\varphi(\frac{|g(u)-g(v)|}{10r_k})).
$$
Hence
\ba
&&|S_{k+1}(Id-S_k)g(w)|\ls\nonumber\\
&&\ls 10 r_k\varphi^{-1}(\frac{r^n}{\eta^{-1}(r_{k+1})^n})
(1+\int_{T}\pint_{B_{k}(u)}\varphi(\frac{|g(u)-g(v)|}{10r_k})\lambda(dv)\lambda(du)).\nonumber
\ea
Let us notice that Lemma \ref{bromba} and the equality $r_{k}=2r_{k+1}$ yield
$$
\frac{1}{8}r_{k}\varphi^{-1}(\frac{r^n}{\eta^{-1}(r_{k+1})^n})\ls\frac{1}{4}r_{k+1}\varphi^{-1}
(\frac{2r^n}{\eta^{-1}(r_{k+1})^n})\ls \ccS_{k+1}.
$$
Take $K_1:=80$. Using Property 1 of $L_{k+1}$, we obtain
\ba
&&|T_{k+1}f(t)-T_k f(t)|
\ls\nonumber\\
&&\ls K_1\ccS_{k+1}(1+\int_{T}\pint_{B_{k}(u)}\varphi(\frac{|R_{k}f(u)-R_{k}f(v)|}{10r_k})
\lambda(dv)\lambda(du)).
\ea
{\bf Case 2.} Fix $k\in J$, $k<m$. It is clear that
$$
T_{k+1}f(t)-T_kf(t)=L_{k+1}(Id-S_{k})S_{k+1} R_{k+1}f(t).
$$
For $k\in J$, the equality $R_k=S_{k+1}R_{k+1}$ holds.
Denote $g:=R_{k}f$, by the definition we get
$$
|(Id-S_{k})g(w)|=|g(w)-\pint_{B_k(w)}g(u)du|.
$$
Property 2 of $R_k$ gives that $g$ is Lipschitz.
Moreover $w\in B_k(w)$, thus due to Corollary \ref{gry0},
we obtain the crucial inequality. For any $A,B>0$
\ba\label{trusia}
&&|g(w)-\pint_{B_k(w)}g(u)du|\ls 6AB
(\int^{\eta^{-1}(r_{k})}_0 \psi(\frac{1}{A\va^{n-1}})\va^{n-1}d\va+\nonumber\\
&&+\frac{1}{n|B_{\|\cdot\|}(0,1)|}\int_{B_{k}(w)}\varphi(\frac{1}{B}\|\nabla g(u)\|_{\ast})du).
\ea
We have used here that $B_{k}(w)=B_{\|\cdot\|}((1-\frac{\eta^{-1}(r_k)}{r})w,\eta^{-1}(r_k))$.
It remains to choose constants $A$ and $B$. For $k\in J$, we have $2\frac{r_k}{\eta^{-1}(r_k)}=
\frac{r_{k+1}}{\eta^{-1}(r_{k+1})}$, so we can take
$$
B:=10\beta \frac{r_{k}}{\eta^{-1}(r_{k})}=5\beta\frac{r_{k+1}}{\eta^{-1}(r_{k+1})},
$$
where $\beta >0$ we choose later. Finding suitable $A$
is more difficult. First we observe that
\ba\label{piol}
&&\int^{\eta^{-1}(r_{k})}_0 \psi(\frac{1}{A\va^{n-1}})\va^{n-1}d\va=
\int^{\eta^{-1}(r_k)}_{2\eta^{-1}(r_{k+1})}\psi(\frac{1}{A\va^{n-1}})\va^{n-1}d\va+\nonumber\\
&&+\int^{2\eta^{-1}(r_{k+1})}_{0}\psi(\frac{1}{A\va^{n-1}})\va^{n-1}d\va.
\ea
Replacing $\va=\frac{\eta^{-1}(r_k)}{r_k}\va'$ 
and applying $2\frac{r_k}{\eta^{-1}(r_k)}=\frac{r_{k+1}}{\eta^{-1}(r_{k+1})}$, we get
$$
\int^{\eta^{-1}(r_k)}_{2\eta^{-1}(r_{k+1})}\psi(\frac{1}{A\va^{n-1}})\va^{n-1}d\va=
\int^{r_{k}}_{r_{k+1}}\psi(\frac{r_k^{n-1}}{A\eta^{-1}(r_k)^{n-1}\va^{n-1}})
\frac{\eta^{-1}(r_k)^n}{r_{k}^{n}}\va^{n-1}d\va.
$$
Since $\eta^{-1}(y)/y$ is increasing,
$\eta^{-1}(\va)^n\ls \frac{\eta^{-1}(r_k)^n}{r_{k}^n}\va^n$, for  $r_{k+1}<\va \ls r_k$.
Due to Lemma \ref{prr} the function $y\psi(1/y)$ is decreasing, so
\be\label{jpw}
\int^{\eta^{-1}(r_k)}_{2\eta^{-1}(r_{k+1})}\psi(\frac{1}{A\va^{n-1}})\va^{n-1}d\va
\ls \int^{r_{k}}_{r_{k+1}}\psi(\frac{\eta^{-1}(r_k)}{r_k}\frac{\va}{A\eta^{-1}(\va)^{n}})
\frac{\eta^{-1}(\va)^n}{\va} d\va.
\ee
We use different estimation for the second integral in (\ref{piol}).
Put $\va=2\eta^{-1}(\va')$, there holds
\ba
&&\int^{2\eta^{-1}(r_{k+1})}_{0}\psi(\frac{1}{A \va^{n-1}})\va^{n-1}d\va=\nonumber\\
&&= 2^n\int^{r_{k+1}}_0\psi(\frac{1}{A2^{n-1}\eta^{-1}(\va)^{n-1}})\eta^{-1}(\va)^{n-1}\eta^{-1}(\va)'d\va.\nonumber
\ea
We can assume that $\eta^{-1}(\va)'$ is right continuous with left limits.
Since $\eta^{-1}(\va)$ is a convex function, the inequality
$\frac{\eta^{-1}(\va)}{\va}\ls \eta^{-1}(\va)'$ holds.
By the convexity of $\psi$ we get
\ba
&&\psi(\frac{1}{A2^{n-1}\eta^{-1}(\va)^{n-1}})\ls
\frac{\eta^{-1}(\va)}{\eta^{-1}(\va)'\va}\psi(\frac{\eta^{-1}(\va)'\va}{A2^{n-1}\eta^{-1}(\va)^{n}})\ls\nonumber\\
&&\ls \frac{\eta^{-1}(\va)}{2^{n-1}\eta^{-1}(\va)'\va}\psi(\frac{\eta^{-1}(\va)'\va}{A\eta^{-1}(\va)^{n}}).\nonumber
\ea
Consequently
\be\label{meta1}
\int^{2\eta^{-1}(r_{k+1})}_{0}\psi(\frac{1}{A\va^{n-1}})\va^{n-1}d\va\ls
2\int^{r_{k+1}}_0\psi(\frac{\eta^{-1}(\va)'\va}{A\eta^{-1}(\va)^{n}})
\frac{\eta^{-1}(\va)^n}{\va} d\va.
\ee
The derivative $\eta^{-1}(\va)'$
can be controlled on the interval $[0,r_{k+1}]$.
Indeed the convexity of $\eta^{-1}$ implies
$$
\eta^{-1}(\va)'\ls \frac{\eta^{-1}(2\va)-\eta^{-1}(\va)}{\va}\ls
\frac{\eta^{-1}(2\va)}{\va}.
$$
Using (\ref{k3}), we obtain
$\eta^{-1}(\va)'\ls 2^{k+2-i}\frac{\eta^{-1}(r_{k})}{r_{k}}$, for $r_{i+1}\ls \va \ls r_{i}$ and $i>k$.
Plugging this estimation into (\ref{meta1}), we obtain
\ba\label{mafia}
&&\int^{2\eta^{-1}(r_{k+1})}_{0}\psi(\frac{1}{A\va^{n-1}})\va^{n-1}d\va\ls\nonumber\\
&&\ls 2\sum^{\infty}_{i=k+1}\int^{r_{i}}_{r_{i+1}}
\psi(2^{k-i}\frac{\eta^{-1}(r_{k})}{r_{k}}\frac{4\va}{A\eta^{-1}(\va)^{n}})
\frac{\eta^{-1}(\va)^n}{\va} d\va.
\ea
We define constant $A$ by the formula
$$
A:=\frac{n\eta^{-1}(r_{k})}{r^n r_{k}}(\ccS_k+4\sum_{i>k}\alpha^{k-i}\ccS_i),
$$
where $\alpha$ must satisfy the condition $1<\alpha<2$.
The definition of $\ccS_k$ and the convexity of $\psi$ give
\ba\label{atena1}
&&\int^{r_{k}}_{r_{k+1}}\psi(\frac{\eta^{-1}(r_k)}{r_k}\frac{\va}{A\eta^{-1}(\va)^{n}})
\frac{\eta^{-1}(\va)^n}{\va} d\va\ls\nonumber\\
&&\ls \int^{r_{k}}_{r_{k+1}}\psi(\frac{r^n \va}{n\ccS_k\eta^{-1}(\va)^{n}})
\frac{\eta^{-1}(\va)^n}{\va} d\va\ls n^{-1}r^n.
\ea
In the same way, for  $i>k$ we prove
\ba\label{atena2}
&&\int^{r_{i}}_{r_{i+1}}
\psi(2^{k-i}\frac{\eta^{-1}(r_{k})}{r_{k}}\frac{4\va}{A\eta^{-1}(\va)^{n}})
\frac{\eta^{-1}(\va)^n}{\va} d\va\ls\nonumber\\
&&
\ls \int^{r_{i}}_{r_{i+1}}
\psi((\frac{\alpha}{2})^{i-k}\frac{r^n\va}{n\ccS_i\eta^{-1}(\va)^{n}})
\frac{\eta^{-1}(\va)^n}{\va} d\va\ls\nonumber\\
&&\ls
n^{-1}(\frac{\alpha}{2})^{i-k}\int^{r_{i}}_{r_{i+1}}
\psi(\frac{r^n\va}{\ccS_i\eta^{-1}(\va)^{n}})
\frac{\eta^{-1}(\va)^n}{\va} d\va=
n^{-1}r^n(\frac{\alpha}{2})^{i-k}.
\ea
Inequalities (\ref{piol}), (\ref{jpw}), (\ref{mafia}), (\ref{atena1}) and (\ref{atena2}) yield
$$
\int^{\eta^{-1}(r_{k})}_0 \psi(\frac{1}{A\va^{n-1}})\va^{n-1}d\va\ls
n^{-1}r^n(1+2\sum_{i>k}(\frac{\alpha}{2})^{i-k})=n^{-1}r^n\frac{2+\alpha}{2-\alpha}.
$$
Denote $\ccS_k':=\ccS_k+4\sum^{\infty}_{i\gs k+1}\alpha^{k-i}\ccS_i$
and $K_2:=60 \beta$.
By (\ref{trusia}) we obtain
\ba\label{haft}
&&|(Id-S_{k})g(w)|\ls\nonumber\\
&& \ls K_2 \ccS_{k}'
( \frac{2+\alpha}{2-\alpha}
+\frac{nr^{-n}}{n|B_{\|\cdot\|}(0,1)|}\int_{B_{k}(w)}
\varphi(\frac{\eta^{-1}(r_{k+1})}{5\beta r_{k+1}}\|\nabla g(u)\|_{\ast} )du)\ls\nonumber\\
&&\ls K_2 \ccS_k'(\frac{2+\alpha}{2-\alpha}+
\pint_{T}
\varphi(\frac{\eta^{-1}(r_{k+1})}{5\beta r_{k+1}}\|\nabla g(u)\|_{\ast} )\lambda(du)),\nonumber
\ea
where we have used the fact that $|B_{\|\cdot\|}(0,1)|=r^{-n}|T|$.
Property 1 of $L_{k+1}$ implies
$$
|T_{k+1}f(t)-T_{k}f(t)|\ls
K_2 \ccS_k'
(\frac{2+\alpha}{2-\alpha}+\pint_{T}
\varphi(\frac{\eta^{-1}(r_{k+1})}{5\beta r_{k+1}}\|\nabla R_{k}f(u)\|_{\ast} )\lambda(du)).
$$
{\bf The last part.} Estimations in cases $k\in I$ and $k\in J$ give
\ba
&&|T_m f(t)-T_0f(t)|\ls\sum^{m-1}_{k=0}|T_{k+1}f(t)-T_k f(t)|\ls \nonumber\\
&&\ls K_1\sum_{k\in I,k< m}\ccS_{k+1}(1+\int_{T}\pint_{B_{k}(u)}\varphi(\frac{|R_{k}f(u)-R_{k}f(v)|}{10r_k})
\lambda(dv)\lambda(du))+\nonumber\\
&&+K_2\sum_{k\in J,k<m}\ccS_k'(\frac{2+\alpha}{2-\alpha}+\pint_{T}
\varphi(\frac{\eta^{-1}(r_{k+1})}{5\beta r_{k+1}}\|\nabla R_{k}f(u)\|_{\ast} )\lambda(du)).\nonumber
\ea
Property 3 of operators $T_k$ gives that $T_0 f$ is a constant function.
Hence
\ba\label{pepe1}
&&|S_m f(s)-S_m f(t)|=|T_m f(s)-T_m f (t)|\ls\nonumber\\
&&\ls 2K_1\sum_{k\in I,k<m}\ccS_{k+1}(1+\int_{T}\pint_{B_{k}(u)}
\varphi(\frac{|R_{k}f(u)-R_{k}f(v)|}{10r_k})
\lambda(dv)\lambda(du))+\nonumber\\
&&+2K_2\sum_{k\in J,k<m}\ccS_k'(\frac{2+\alpha}{2-\alpha}+\pint_{T}
\varphi(\frac{\eta^{-1}(r_{k+1})}{5\beta r_{k+1}}\|\nabla R_{k}f(u)\|_{\ast} )\lambda(du)).
\ea
Let us notice that $B_k(u)\subset B(u,2r_k)$. We use the above inequality to
prove that for each process $X(t)$, $t\in T$ which satisfies (\ref{war2}) the inequality
holds
$$
\E\sup_{s,t\in T}|X(s)-X(t)|\ls K\sum_{k\gs 0}\ccS_k,
$$
where the constant $K$ depends only on $n$.
Due to the remark (\ref{hlyt}) and Lemma \ref{jujki}
we can assume that a process $X(t)$, $t\in T$
has finite number of different Lipschitz trajectories.
\bl\label{adas}
Let $0\ls k< m$. For each
$u,v\in T$, where  $d(u,v)\ls 2r_k$, we have
$$
\E\varphi(\frac{|R_k X(u)- R_k X(v)|}{10r_k})\ls 1.
$$
\el
\begin{dwd}
We denote $X_k:=R_k X$.
If $R_{k}=Id$, then the condition (\ref{war2}) implies the lemma. Otherwise there exists
$N>0$ and a sequence $k=k_0<k_1<...<k_N\ls m$
such that $X_{k}=R_kX=S_{k_1}S_{k_2}...S_{k_N}X$.
For  simplicity we denote $u_{k_0}:=u$, $v_{k_0}:=v$.
We obtain the following equalities:
\ba
&& X_k(u)=\pint_{B_{k_1}(u_{k_0})}\pint_{B_{k_2}(u_{k_1})}...\pint_{B_{k_N}(u_{k_{N-1}})}X(u_{k_N})
\lambda(du_{k_N})...\lambda(du_{k_1}); \nonumber\\
&& \label{klapa2}
X_k(v)=\pint_{B_{k_1}(u_{k_0})}\pint_{B_{k_2}(v_{k_1})}...\pint_{B_{k_N}(v_{k_{N-1}})}X(v_{k_N})
\lambda(dv_{k_N})...\lambda(dv_{k_1}).\nonumber
\ea
Take $u_{k_{i+1}}\in B_{k_i}(u_{k_i})\subset B(u_{k_i},2r_{k_i})$.
By the triangle inequality and (\ref{k3}) we have
$$
d(u_{k_0},u_{k_N})\ls\sum^{N-1}_{i=0}d(u_{k_i},u_{k_{i+1}})\ls 2\sum^{N-1}_{i=0}r_{k_i}\ls 4r_{k_0},
\;\;d(v_{k_0},v_{k_N})\ls 4r_{k_0}.
$$
That means, for some probability measures
$\nu_{u},\nu_{v}$ with supports respectively in $B(u,4r_k)$, $B(v,4r_k)$ the equalities hold
$$
X_k(u)=\int_{B(u,4r_k)}X(w)\nu_{u}(dw),\;\;
X_k(v)=\int_{B(v,4r_k)}X(z)\nu_{v}(dz).
$$
By the assumption $d(u,v)\ls 2r_k$, hence $d(w,z)\ls 10r_k$.
The Jensen inequality, the Fubini theorem and (\ref{war2}) yield
\ba
&&\E\varphi(\frac{|X_k(u)-X_k(v)|}{10r_k})\ls\nonumber\\
&&\ls \pint_{B(u,4r_k)}\pint_{B(v,4r_k)}\E\varphi(\frac{|X(w)-X(z)|}{d(w,z)})
\nu_{u}(dw)\nu_{v}(dz)\ls 1.\nonumber
\ea
\end{dwd}
The Auerbach lemma (for the proof see \cite{Woj} - Lemma 11, II.E.)
gives that there exists biorthogonal system $((b_i,b_i^{\ast}))$
in the space $\R^{n}\times \R^{n}$ such that $\|b_i\|=1,\|b^{\ast}_i\|_{\ast}=1$.
Consequently for each $v\in \R^n$
\be\label{vapi}
\|v\|_{\ast}=
\sup_{u\in B_{\|\cdot\|}(0,1)}|\langle v, u\rangle|\ls
\sup_{u\in B_{\|\cdot\|}(0,1)}\sum^n_{i=1}|\langle v, b_i \rangle||\langle b_i^{\ast},u\rangle|\ls
\sum^n_{i=1}|\langle v, b_i \rangle|.
\ee
\bl\label{gry}
For $|\cdot|$-almost all $t\in T$ the following inequality holds
$$
\E\varphi(\frac{\eta^{-1}(r_{k+1})}{5\beta r_{k+1}}\|\nabla R_k X(t)\|_{\ast})\ls 1,\;\;0\ls k<m, k\in J,
$$
where $\beta:=\sum^n_{i=1}\beta(b_i)$ ($\beta_i$ was defined in Lemma \ref{bagi}). 
Since $\beta_i\ls n$, thus $\beta\ls n^2$.
\el
\begin{dwd}
We put $X_k:=R_{k} X$.
The definition gives that $X_k=S_{k+1}R_{k+1}X=(X_{k+1})_{r_{k+1}}$, for $k\in J$.
Applying (\ref{vapi}) and Lemma \ref{bagi} we obtain that for $|\cdot|$-almost all $t\in T$
there holds  
\ba
&&\frac{\eta^{-1}(r_{k+1})}{5\beta r_{k+1}}\|\nabla X_k (t)\|_{\ast}\ls
\sum^n_{i=1}\frac{1}{\beta}\frac{\eta^{-1}(r_{k+1})}{5r_{k+1}}|\langle \nabla  (X_{k+1})_{r_{k+1}}(t), b_i\rangle |\ls   \nonumber\\
&&\ls \sum^n_{i=1}\frac{\beta(b_i)}{\beta}
\pint_{\triangle^{b_i}_{t,r_{k+1}}}\pint_{\triangle^{-b_i}_{t,r_{k+1}}}\frac{|X_{k+1}(u)-X_{k+1}(v)|}
{10r_{k+1}}\sigma^i_1(du)
\sigma^i_2(dv),\nonumber
\ea
where $\sigma^i_1(du)=\sigma^{b_i}_{t,r_{k+1}}(du)$, $\sigma^i_2(dv)=\sigma^{-b_i}_{t,r_{k+1}}(dv)$.
The Jensen inequality yields
\ba
&&\varphi(\frac{\eta^{-1}(r_{k+1})}{5\beta r_{k+1}}\|\nabla X_{k}(t)\|_{\ast})\ls\nonumber\\
&&\ls \sum^n_{i=1}\frac{\beta(b_i)}{\beta}
\pint_{\triangle^{b_i}_{t,r_{k+1}}}\pint_{\triangle^{-b_i}_{t,r_{k+1}}}
\varphi(\frac{|X_{k+1}(u)-X_{k+1}(v)|}{10r_{k+1}})
\sigma^i_1(du)
\sigma^i_2(dv).\nonumber
\ea
Notice that $d(u,v)\ls 2r_{k+1}$ for $u\in \triangle^{b_i}_{t,r_{k+1}}$,
$v\in \triangle^{-b_i}_{t,r_{k+1}}$. The Fubini theorem and Lemma \ref{adas} imply
$$
\E\varphi(\frac{\eta^{-1}(r_{k+1})}{5\beta r_{k+1}}
\|\nabla X_{k}(t)\|_{\ast})\ls \sum^n_{i=1}\frac{\beta(b_i)}{\beta}=1.
$$
\end{dwd}
By the Fubini theorem, Lemma \ref{adas}, Lemma \ref{gry} and (\ref{pepe1}) we obtain
$$
\E\sup_{s,t\in T}|S_m X(s)-S_m X(t)|\ls 4K_1\sum_{k\in I,k<m}\ccS_{k+1}+
2K_2(1+\frac{2+\alpha}{2-\alpha})\sum_{k\in J,k<m}\ccS_k'.
$$
Let us notice that $\sum_{k\in J,k<N}\ccS_k'\ls (1+\frac{4}{\alpha-1})\sum_{k\gs 0}\ccS_k$.
We have proved that for some constant $K$ which depends only $n$ the following inequality holds
$$
\E\sup_{s,t\in T}|S_m X(s)-S_m X(t)|\ls K\sum_{k\gs 0}\ccS_k.
$$
Since $\lim_{m\ra \infty}S_m X(t)=X(t)$, thus due to the Fatou lemma
$$
\E\sup_{s,t\in T}|X(s)-X(t)|\ls
\liminf_{m\ra \infty}\E \sup_{s,t\in T}|S_m X(s)-S_m X(t)|\ls
K\sum_{k\gs 0}\ccS_k.
$$
It ends the proof of the theorem.
\end{dw2}

\section{The case of $\eta'(0)<\infty$}

We assume that $\eta'(0)<\infty$.
Let us remind that there exists $m\gs 0$ such that
$r_m>0$ and $r_{m+1}=0$. We have defined $\ccS_m$ as the infimum over all $c>0$ such that
$$
\int^{r_m}_{0}\frac{\lambda(B(0,\va))}{\va}\psi(\frac{\va}{c\lambda(B(0,\va))})d\va=
\int^{r_m}_{0}\frac{\eta^{-1}(\va)^n}{r^n\va}\psi(\frac{r^n\va}{c \eta^{-1}(\va)^n})d\va\ls 1
$$
and $\ccS_{k}=0$, for $k>m$.
\begin{dw3}
We follow the proof of Theorem \ref{teo1} in the case of $\eta'(0)=\infty$. 
The only difference is when $\ccS_m=\infty$. For all $0<\delta\ls r_m/2$
we denote by $\ccS_m(\delta)$ numbers such that 
$$ 
\int^{r_m}_{\delta}\frac{\lambda(B(0,\va))}{\va}\psi(\frac{\va}{\ccS_m^\delta \lambda(B(0,\va))})d\va=
\int^{r_m}_{\delta}\frac{\eta^{-1}(\va)^n}{r^n\va}\psi(\frac{r^n\va}{\ccS_m^\delta \eta^{-1}(\va)^n})d\va=1.
$$
Since $\ccS_m=\infty$ we have $\lim_{\delta\ra\infty} \ccS_{m}(\delta)=\infty$. The proof
of the left-hand side of (\ref{kwiat1}) in Theorem \ref{teo1} implies
$\ccS_m(\delta)\ls 3(n+2)S(T,d,\varphi)$. Hence $S(T,d,\varphi)=\infty$. 
It ends the proof. 
\end{dw3}
\smallskip

\noindent
{\bf Acknowledgment} I would like to thank professor Stanislaw Kwapien and 
the anonymous referee for numerous remarks which
helped me to improve the paper.

\flushleft{\textsc{W. Bednorz\\
Department of Mathematics\\
University of Warsaw\\
Banacha 2, 02-097 Warsaw, Poland}\\
\textit{E-mail address: wbednorz@mimuw.edu.pl}}

\end{document}